\documentclass[11pt]{article}
\usepackage[utf8]{inputenc}

\usepackage[letterpaper,margin=1in]{geometry}

\usepackage{graphicx}

\usepackage{xcolor}
\definecolor{darkblue}{rgb}{0,0,0.5}
\definecolor{darkred}{rgb}{.8,0,0}
\definecolor{darkgreen}{rgb}{0,.5,0}
\usepackage[colorlinks=true,linkcolor=darkblue,citecolor=darkblue,urlcolor=darkblue]{hyperref}

\usepackage{amsthm, amsmath, amssymb} 
\usepackage{enumerate} 

\usepackage{amsfonts}
\usepackage{graphicx}
\usepackage{placeins}
\usepackage{floatrow}
\usepackage{lipsum} 

\usepackage{titlesec}
\usepackage{fancyhdr}

\usepackage{tikz}
\usetikzlibrary{shapes}
\usetikzlibrary{arrows}
\usetikzlibrary{calc}
\usepackage{subcaption}

\usepackage{tikz}
\usetikzlibrary{shapes, arrows, calc, positioning}

\usepackage[shortlabels]{enumitem} 

\usepackage{multicol}
\usepackage[mathscr]{eucal}
\usepackage{subcaption}

\usepackage{color}
\usepackage{marginnote}

\usepackage{graphicx}

\makeatletter
\newtheorem*{rep@theorem}{\rep@title}
\newcommand{\newreptheorem}[2]{%
\newenvironment{rep#1}[1]{%
 \def\rep@title{#2 \ref{##1}}%
 \begin{rep@theorem}}%
 {\end{rep@theorem}}}
\makeatother

\theoremstyle{plain}
\newtheorem{theorem}{Theorem}
\newreptheorem{theorem}{Theorem} 
\newtheorem*{theorem*}{Theorem}
\newtheorem{lemma}[theorem]{Lemma}
\newtheorem*{lemma*}{Lemma}

\newtheorem*{proposition*}{Proposition}

\newtheorem*{corollary*}{Corollary}

\newtheorem*{claim*}{Claim}
\newtheorem{conjecture}{Conjecture}
\newtheorem*{conjecture*}{Conjecture}

\theoremstyle{definition}
\newtheorem{definition}{Definition}
\newtheorem*{definition*}{Definition}

\theoremstyle{remark}

\newtheorem*{example*}{Example}

\newtheorem{note}{Note}

\newcommand{\floor}[1]{\left\lfloor #1 \right\rfloor}

\renewcommand{\natural}{\mathbb{N}} 

\renewcommand{\bar}{\overline}

\renewcommand{\bar}{\overline}

\newcommand{\FC}{\mathfrak{C}}

\newcommand{\roottwo}{1.4142}

\newcommand{\tikzArrow}{\draw[-latex] (ArrowHead) -- (ArrowTail)}
\newcommand{\vertremoved}{black}
\newcommand{\recolored}{gray}
\newcommand{\edgeremoved}{dashed}

\newcommand{\GsquaredChoose}{Therefore, $G^2$ can be colored from any lists of size 12, which contradicts the assumption that $\chi_\ell(G^2) > 12$. }

\newcommand{\LemmaReduce}{We can see that the conditions of Lemma \ref{lem:reduceWegner4} are satisfied. }

\tikzstyle{vertex}=[draw,thick,fill=white,circle,inner sep=2pt]

\def\blob#1#2#3{\draw[fill=#3,rounded corners=#1*3mm] (#2) +($(0:#1*2+#1*rnd)$)
\foreach \a in {20,40,...,350} {  -- +($(\a: #1*2+#1*rnd)$) } -- cycle;}

\title{Square of Planar Graphs of Max Degree Four without Five Cycles}
\author{Eric Culver\\
Department of Mathematical and Statistical Sciences\\
University of Colorado Denver\\
\texttt{eric.culver@ucdenver.edu}
\and
Stephen G. Hartke\thanks{Supported in part by a Collaboration Grant from the Simons Foundation (\#316262 to Stephen G. Hartke).}\\
Department of Mathematical and Statistical Sciences\\
University of Colorado Denver\\
\texttt{stephen.hartke@ucdenver.edu}
}
\date{May 2022}

\begin{document}

\maketitle

\begin{abstract}
    We show that the choosability of the square of planar graphs of max degree 4 without five cycles is at most 12.
    
    
\end{abstract}

\section{Introduction}

In 1977, Wegner conjectured the following upper bounds on the chromatic numbers of squares of planar graphs:

\begin{conjecture}[Wegner~\cite{Wegner1977preprint}]
    Let $G$ be a planar graph with maximum degree $\Delta$. Then
    \[ \chi(G^2) \le \begin{cases} 7 & \Delta \le 3 \\ \Delta + 5 & 4 \le \Delta \le 7 \\ \floor{\frac{3\Delta}{2}} + 1 & \Delta \ge 8 \end{cases} \]
\end{conjecture}

In this paper, we will be focusing on the $\Delta = 4$ case.

\begin{conjecture}
    Let $G$ be a planar graph with maximum degree 4, then $\chi(G^2) \le 9$.
\end{conjecture}

\begin{figure}[ht]
    \centering
    \begin{tikzpicture}[node distance=2cm, every node/.style={draw=black, circle}]
        \node (A) at (0,0) {};
        \node (B) at (1,0) {};
        \node (C) at (0,1) {};
        \node (D) at (-1,0) {};
        \node (E) at (0,-1) {};
        \node (F) at (2,2) {};
        \node (G) at (-2,2) {};
        \node (H) at (2,-2) {};
        \node (I) at (-2,-2) {};
        \draw (B) -- (A) -- (D);
        \draw (C) -- (A) -- (E);
        \draw (F) -- (B) -- (H);
        \draw (F) -- (C) -- (G);
        \draw (G) -- (D) -- (I);
        \draw (H) -- (E) -- (I);
        \draw (F) -- (G) -- (I) -- (H) -- (F);
    \end{tikzpicture}
    \caption{Sharpness Example}
    \label{fig:sharpness}
\end{figure}
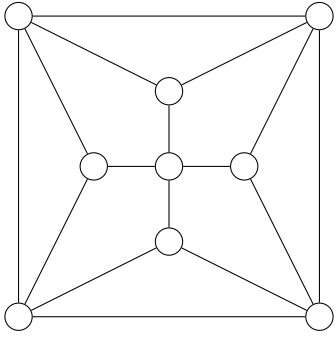

We can see that this conjecture, if true, would be sharp, from the graph in Figure \ref{fig:sharpness}, which is a planar graph of maximum degree four on nine vertices. We can see that the square of this graph is the complete graph on nine vertices. Also, since this graph contains vertices of degree three, there is an infinite number of examples of planar graphs of maximum degree four which contain this graph as a subgraph, and therefore would also require nine colors in their squares.

Not much progress was made on the $\Delta = 4$ case of this conjecture until 2002, when Borodin et al.~\cite{Borodin2002} showed that for planar graphs $G$ with maximum degree $\Delta \le 20$ that $\chi(G^2) \le 59$. This was improved by Zhu and Bu~\cite{Zhu2018} in 2018 who showed that planar graphs $G$ with maximum degree $\Delta \le 5$ have $\chi(G^2) \le 20$.

In this paper, we will be looking at the more general case of \emph{list coloring}, which was first introduced independently by Vizing in 1976~\cite{Vizing1976} and by Erd\H{o}s, Rubin, and Taylor in 1979~\cite{ErdosRubinTaylor1979}. A \emph{list assignment} for a graph $G$ is a function $L$ that assigns to each vertex a list of colors. An \emph{$L$-coloring} of $G$ is a coloring of $G$ such that for each vertex, the color assigned to it is picked from its list. A graph $G$ is \emph{$k$-choosable} if there exists an $L$-coloring of $G$ for every assignment $L$ of lists of size $k$ to the vertices of $G$.

\begin{figure}[ht]
    \centering
    \begin{tikzpicture}[node distance=2cm, every node/.style={draw=black, circle}]
        \node (A) at (0,1.5) {};
        \node (B) at (0,0.5) {};
        \node (C) at (0,-0.5) {};
        \node (D) at (0,-1.5) {};
        \node (E) at (2,0) {};
        \node (F) at (-2,0) {};
        \draw (E) -- (A) -- (F);
        \draw (E) -- (B) -- (F);
        \draw (E) -- (C) -- (F);
        \draw (E) -- (D) -- (F);
    \end{tikzpicture}
    \caption{2-colorable, but not 2-choosable}
    \label{fig:wegner_K_24}
\end{figure}
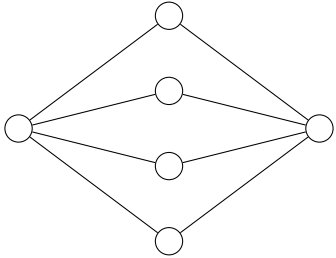

Every $k$-choosable graph is $k$-colorable. However, the converse is known to not be true in general. For example, the graph in Figure \ref{fig:wegner_K_24} is 2-colorable but not 2-choosable. 

The \emph{choosability}, $\chi_\ell(G)$, of $G$ is the minimum $k$ such that $G$ is $k$-choosable. Then, this relationship is expressed as:
\[ \chi(G) \le \chi_\ell(G) \text{ for all graphs $G$} \]
Where the inequality is sometimes strict.

In this paper, we will be proving Theorem \ref{thm:wegner4no5}.

\begin{reptheorem}{thm:wegner4no5}
	Let $G$ be a planar graph with no 5-cycles such that $\Delta(G) \le 4$. Then $\chi_\ell(G^2) \le 12$.
\end{reptheorem}

Our technique for proving this statement will be the discharging method. The method of discharging was developed by Birkoff~\cite{Birkhoff1913} and Heesch~\cite{Heesch1969, Heesch1972}, and was ultimately used by Appel, Haken, and Koch~\cite{AppelHaken1977, AppelHakenKoch1977b, AppelHaken1989} to prove the Four Color Theorem. It has since been used to great effect for many results in graph theory.

For a survey of results proven by discharging, I refer the reader to~\cite{Borodin2013, JendrolVoss2013}. For a more thorough explanation of the discharging method than this paper will provide, and a further survey of results, see~\cite{cranston_west_2017}.

All discharging arguments follow the same steps:
\begin{enumerate}[a.]
    \item We suppose our graph $G$ is a minimal counterexample to the statement
    we want to prove. \label{discharge:minCE_color}
    \item We argue that $G$ cannot contain certain configurations, called \emph{reducible} configurations. If $G$ did contain a reducible configuration, then we can construct a smaller graph $G'$ for which the statement holds by the minimality of $G$. We then show that the statement (usually a coloring) can be extended from $G'$ to $G$. This then shows the statement holds for $G$, which is a contradiction. Therefore, $G$ cannot contain any of these reducible configurations. \label{discharge:reducibility_color}
    \item We then use the technique of discharging to show that $G$ must contain one of the reducible configurations. We then call these configurations \emph{unavoidable}. \label{discharge:discharging_color}
    \item This is a contradiction. Therefore, the statement is true.
\end{enumerate}

Our proof will follow this same basic outline. 
We will handle Step \ref{discharge:reducibility_color} in Section \ref{sec:color_reduce}, where our chief tool will be Lemma \ref{lem:reduceWegner4}.
We will handle Step \ref{discharge:discharging_color} in Section \ref{sec:color_discharge}. 

In order to show this, we will need some definitions.

Let $[k] = \{0, 1, \cdots, k-1 \}$.

\begin{definition}
In a graph $G$, a \emph{$k$-vertex} is a vertex of degree $k$. Similarly, a \emph{$k^+$-vertex} is a vertex of degree at least $k$, and a \emph{$k^-$-vertex} is a vertex of degree at most $k$.

A \emph{$k$-face}, \emph{$k^+$-face}, and \emph{$k^-$-face} are defined similarly, referring to the length of the face instead of the degree of the vertex.
\end{definition}
\begin{definition}
Given a set $S \subseteq V(G)$ of vertices of $G$, the \emph{neighborhood} $N(S)$ of $S$ is the set of vertices in $G \setminus S$ which are adjacent to at least one vertex in $S$. 

The \emph{closed neighborhood} $\bar{N}(S)$ of $S$ includes $S$, i.e., $\bar{N}(S) = N(S) \cup S$.

If the set $S$ is small, we often omit the curly braces, so: $N(x) = N(\{x\})$ and $N(x,y) = N(\{x,y\})$.
\end{definition}

\begin{definition}
	Given a plane graph $G$, the set of faces of $G$ will be notated $F(G)$.
\end{definition}

\section{Reducibility}
\label{sec:color_reduce}

The intuitive idea behind the following lemma is that $X$ is the set of vertices and $Y$ the set of edges that are removed from $G$ to produce $H$. The vertices in $R$ are being recolored. The first condition of the lemma simply says that if we are removing a vertex, we also must remove all the edges incident to that vertex. The second condition ensures that the coloring on $H^2$ is a valid coloring of the vertices of $P$ in $G^2$. The third condition then checks that we can extend the coloring on the vertices in $P$ to the vertices in $X$ and $R$.  

\begin{lemma}
	\label{lem:reduceWegner4}
	Let $G$ be a graph with vertex set $V(G)$ and edge set $E(G)$.
	Let $X, R, P \subseteq V(G)$ be three disjoint subsets of the vertex set of $G$ whose union is $V(G)$.
	And let $Y, Q \subseteq E(G)$ be two disjoint subsets of the edge set of $G$ whose union is $E(G)$.
	Let $H$ be the subgraph of $G$ on vertex set $R \cup P$ and edge set $Q$.
	If $X, R, P, Y, Q, H$, and $G$ satisfy:
	\begin{enumerate}
		\item Any edge of $G$ incident to a vertex in $X$ must be in $Y$.
		\item Any edges in $G^2$ not in $H^2$ must be incident to vertices in $X \cup R$.
		\item Define $f: X \cup R \to \natural$ by $f(v) = 12 - |N_{G^2}(v) \cap P|$. Then the subgraph of $G^2$ induced by $X \cup R$ is $f$-choosable.
		\item $\chi_\ell(H^2) \le 12$
	\end{enumerate}
	Then $\chi_\ell(G^2) \le 12$.
\end{lemma}
\begin{proof}
	Let $L$ be a list assignment for $G$, mapping each vertex to a list of 12 colors. Since $H$ is a subgraph of $G$, this is also a list assignment for $H$. Since $\chi_\ell(H^2) \le 12$, $H$ is $L$-square-colorable. Since the only edges of $G^2$ not in $H^2$ are incident to vertices in $X \cup R$, if we erase the color on vertices in $R$, then extend the coloring to the vertices in $X \cup R$, we will have square colored $G$. The condition that the subgraph of $G^2$ on vertex set $X \cup R$ is $f$-choosable for those specific $f$-values is sufficient to extend this coloring. Therefore, $G$ is $L$-square-colorable.
	Since $L$ was arbitrarily chosen, this shows $\chi_\ell(G^2) \le 12$.
\end{proof}

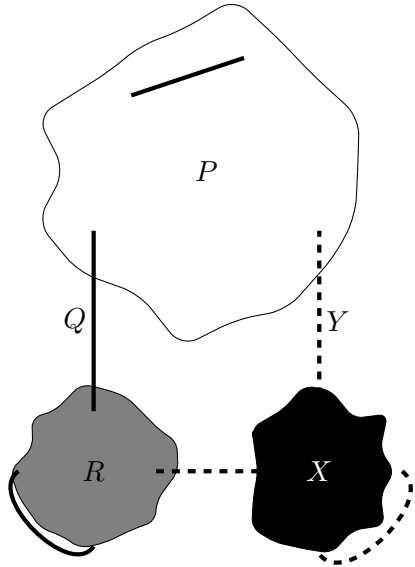
\begin{figure}[ht]
    \begin{center}
        \begin{tikzpicture}
            \blob{0.4}{0,0}{black}
            \blob{0.4}{-3,0}{gray}
            \blob{0.8}{-1.5,4}{white}
            \node[text=white] at (0,0) {$X$};
            \node at (-3,0) {$R$};
            \node at (-1.5,4) {$P$};
            \draw[dashed, ultra thick] (0,0.8) -- (0,3.2);
            \draw[dashed, ultra thick] (1.1,0) edge[bend left=100] (0,-1.1);
            \draw[dashed, ultra thick] (-0.8,0) -- (-2.2,0);
            \draw[ultra thick] (-3,0.8) -- (-3,3.2);
            \draw[ultra thick] (-4,0) edge[bend right=100] (-3,-1);
            \draw[ultra thick] (-2.5,5) -- (-1,5.5);
            \node at (0.25,2) {$Y$};
            \node at (-3.25,2) {$Q$};
        \end{tikzpicture}
        \caption{Guide for Lemma \ref{lem:reduceWegner4}}
        \label{fig:reduceWegner4}
    \end{center}
\end{figure}

Note that in all of our applications of this lemma except one, we show that all the elements of $X \cup R$ must be distinct, and the subgraph of $G^2$ on $X \cup R$ is complete. This means we need only consider the configurations in the most general way possible, since any amount of overlap in the configuration will only increase the $f$-values, and it cannot add any edges into the subgraph on $X \cup R$, since that is already complete.

\begin{theorem}
	\label{thm:wegner4no5}
	Let $G$ be a planar graph with no 5-cycles such that $\Delta(G) \le 4$. Then $\chi_\ell(G^2) \le 12$.
\end{theorem}

Let $\FC$ be the family of all planar graphs with no 5-cycles such that $\Delta(G) \le 4$.
For the rest of this paper, let $G$ be a minimal counterexample to Theorem \ref{thm:wegner4no5}, that is, $G \in \FC$ and $\chi_\ell(G^2) > 12$ and for any graph $H \in \FC$ with a smaller number of edges than $G$ or a smaller number of vertices than $G$, $\chi_\ell(H^2) \le 12$. 

From these assumptions, we can prove certain lemmas about $G$.
All of these lemmas will take the form: ``The graph $G$ cannot contain structure $\mathcal{X}$'', and most will be proven by assuming that $G$ does have structure $\mathcal{X}$ and finding $X, R, P, Y, Q, H$ which satisfy the conditions of Lemma \ref{lem:reduceWegner4}. 
Since this lemma concludes that $\chi_\ell(G^2) \le 12$, while by assumption $\chi_\ell(G^2) > 12$, this leads to a contradiction, showing that $G$ cannot have structure $\mathcal{X}$. 
Note that to derive the condition that $\chi_\ell(H^2) \le 12$ we will use that $H \in \FC$ and that $H$ is a smaller graph than $G$, meaning that at least one of $X, Y$ must be nonempty. These are further conditions we will need to check. 

\begin{note}
	In the following figures, we will mark the removed vertices (elements of $X$) by filling them in with black, the removed edges (elements of $Y$) by dashing them, and recolored vertices (element of $R$) by filling them in with gray. We need only show the unremoved edges and precolored vertices of $G$ that are within distance two of an element of $X$, $Y$, or $R$, as those are the only ones that contribute to the count of the $f$-values. We will be showing those in the most general way possible, as if nothing overlapped, and every vertex was of maximum degree. To facilitate the checking of Lemma \ref{lem:reduceWegner4} by the reader, the subgraph of $G^2$ on the vertex set $X \cup R$ with the $f$ values is also given, after an arrow.
\end{note}

\begin{lemma}
    \label{lem:conn}
    The graph $G$ is connected.
\end{lemma}
\begin{proof}
    Suppose $G$ did have multiple connected components, and let $G_1$ be one of them, while $G_2$ is the rest of the graph. 
    By the minimality of $G$, we can assume that $\chi_\ell(G_1^2) \le 12$ and $\chi_\ell(G_2^2) \le 12$. Since they are disconnected, given any assignment of lists of 12 colors to the vertices of $G$, we can color $G_1$, without affecting $G_2$, and then color $G_2$, without affecting $G_1$. In this way, we can always color $G$. 
    \GsquaredChoose
    Therefore, $G$ is connected.
\end{proof}

\begin{figure}[ht]
    \begin{center}
    \begin{tikzpicture}[node distance=2cm, every node/.style={draw=black, circle}]
        \coordinate[draw=none] (GraphShift) at (7,0) {};
        \coordinate[draw=none] (ArrowHead) at (4,0) {};
        \coordinate[draw=none] (ArrowTail) at (5,0) {};
        \coordinate[draw=none] (NumShift) at (0.35,0.35) {};
        \node[fill=\vertremoved] (A) at (0:1) {};
        \node (B) at (1*360/2:1) {};
        \node (Bm) at ($(B) + (90:1)$) {};
        \node (Bn) at ($(B) + (180:1)$) {};
        \node (Bo) at ($(B) + (270:1)$) {};
        \node[fill=\vertremoved] (Az) at ($(A) + (GraphShift)$) {};
        \node[draw=none] (AzLabel) at ($(Az) + (NumShift)$) {8};
        \draw[\edgeremoved] (A) -- (B);
        \draw (Bm) -- (B) -- (Bn) (Bo) -- (B);
        \tikzArrow;
    \end{tikzpicture}
    \caption{A 1-vertex}
    \label{fig:no1v}
    \end{center}
\end{figure}
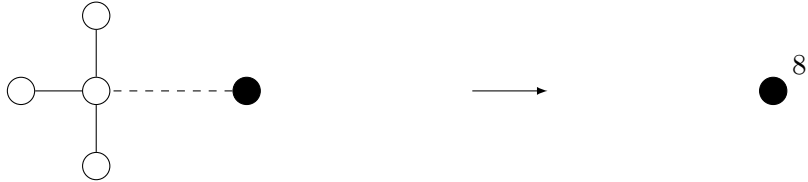

\begin{lemma}
    \label{lem:no1v}
    The graph $G$ cannot have a 1-vertex.
\end{lemma}
\begin{proof}
    Suppose it did have a 1-vertex $v$.
	Let $X = \{v\}$, let $Y$ be the edge incident to $v$, and let $R$ be empty. (See Figure \ref{fig:no1v})
	\LemmaReduce
	\GsquaredChoose
	Therefore, $G$ cannot have a 1-vertex.
\end{proof}

\begin{figure}[ht]
    \begin{center}
    \begin{tikzpicture}[node distance=2cm, every node/.style={draw=black, circle}]
        \node[fill=\vertremoved] (A) at (0:1) {};
        \node (B) at (1*360/3:1) {};
        \node (C) at (2*360/3:1) {};
        \node (Bm) at ($(B) + (90:1)$) {};
        \node (Bn) at ($(B) + (150:1)$) {};
        \node (Cm) at ($(C) + (-90:1)$) {};
        \node (Cn) at ($(C) + (-150:1)$) {};
        \node[fill=\vertremoved] (Az) at ($(A) + (GraphShift)$) {};
        \node[draw=none] (AzLabel) at ($(Az) + (NumShift)$) {6};
        \draw[\edgeremoved] (C) -- (A) -- (B);
        \draw (B) -- (C);
        \draw (Bm) -- (B) -- (Bn);
        \draw (Cm) -- (C) -- (Cn);
        \tikzArrow;
    \end{tikzpicture}
    \caption{A 2-vertex incident to a 3-face}
    \label{fig:no2v3f}
    \end{center}
\end{figure}
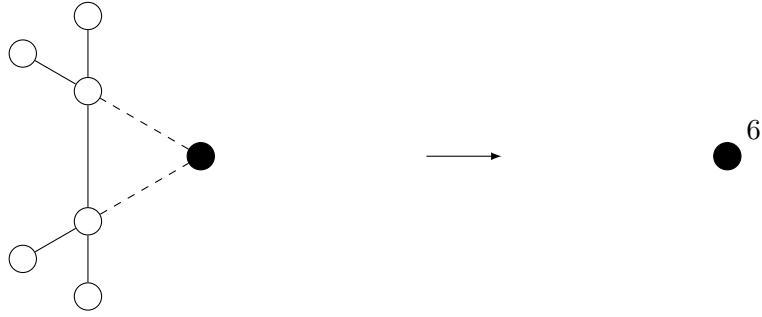

\begin{lemma}
	\label{lem:no2v3f}
	The graph $G$ cannot have a 2-vertex incident to a 3-face.
\end{lemma}
\begin{proof}
	Suppose it did have a 2-vertex $v$ incident to a 3-face.
	Let $X = \{v\}$, let $Y$ be the two edges incident to $v$, and let $R$ be empty. (See Figure \ref{fig:no2v3f})
	\LemmaReduce
	\GsquaredChoose
	Therefore, $G$ cannot have a 2-vertex incident to a 3-face.
\end{proof}

\begin{figure}[ht]
    \begin{center}
    \begin{tikzpicture}[node distance=2cm, every node/.style={draw=black, circle}]
        \node[fill=\vertremoved] (A) at (0:1) {};
        \node (B) at (1*360/4:1) {};
        \node (C) at (2*360/4:1) {};
        \node (D) at (3*360/4:1) {};
        \node (Bm) at ($(B) + (45:1)$) {};
        \node (Bn) at ($(B) + (135:1)$) {};
        \node (Dn) at ($(D) + (-45:1)$) {};
        \node (Dm) at ($(D) + (-135:1)$) {};
        \node[fill=\vertremoved] (Az) at ($(A) + (GraphShift)$) {};
        \node[draw=none] (AzLabel) at ($(Az) + (NumShift)$) {5};
        \draw[\edgeremoved] (D) -- (A) -- (B);
        \draw (B) -- (C) -- (D);
        \draw (Bm) -- (B) -- (Bn);
        \draw (Dm) -- (D) -- (Dn);
        \tikzArrow;
    \end{tikzpicture}
    \caption{A 2-vertex incident to a 4-face}
    \label{fig:no2v4f}
    \end{center}
\end{figure}
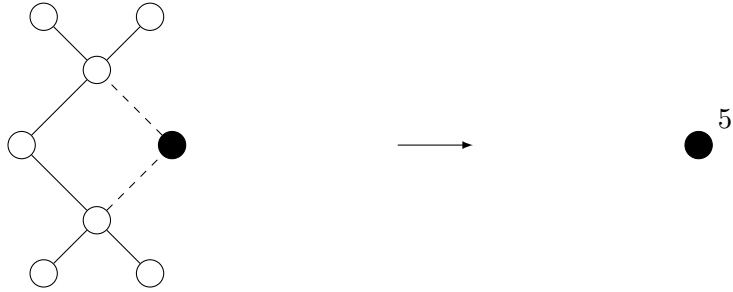

\begin{lemma}
	\label{lem:no2v4f}
	The graph $G$ cannot have a 2-vertex incident to a 4-face.
\end{lemma}
\begin{proof}
	Suppose it did have a 2-vertex $v$ incident to a 4-face.
	Let $X = \{v\}$, let $Y$ be the two edges incident to $v$, and let $R$ be empty. (See Figure \ref{fig:no2v4f})
	\LemmaReduce
	\GsquaredChoose
	Therefore, $G$ cannot have a 2-vertex incident to a 4-face.
\end{proof}

\begin{figure}[ht]
    \begin{center}
    \begin{tikzpicture}[node distance=2cm, every node/.style={draw=black, circle}]
        \node[fill=\vertremoved] (A) at (0*360/2:1) {};
        \node[fill=\vertremoved] (B) at (1*360/2:1) {};
        \node (Am) at ($(A) + (0:1)$) {};
        \node (Amm) at ($(Am) + (30:1)$) {};
        \node (Amn) at ($(Am) + (0:1)$) {};
        \node (Amo) at ($(Am) + (-30:1)$) {};
        \node (Bm) at ($(B) + (180:1)$) {};
        \node (Bmm) at ($(Bm) + (180+30:1)$) {};
        \node (Bmn) at ($(Bm) + (180+0:1)$) {};
        \node (Bmo) at ($(Bm) + (180-30:1)$) {};
        \node[fill=\vertremoved] (Az) at ($(A) + (GraphShift)$) {};
        \node[draw=none] (AzLabel) at ($(Az) + (NumShift)$) {7};
        \node[fill=\vertremoved] (Bz) at ($(B) + (GraphShift)$) {};
        \node[draw=none] (BzLabel) at ($(Bz) + (NumShift)$) {7};
        \draw[\edgeremoved] (A) -- (B);
        \draw[\edgeremoved] (Am) -- (A);
        \draw (Amm) -- (Am) -- (Amn) (Am) -- (Amo);
        \draw[\edgeremoved] (Bm) -- (B);
        \draw (Bmm) -- (Bm) -- (Bmn) (Bm) -- (Bmo);
        \draw (Az) -- (Bz);
        \tikzArrow;
    \end{tikzpicture}
    \caption{Adjacent 2-vertices}
    \label{fig:no22v}
    \end{center}
\end{figure}
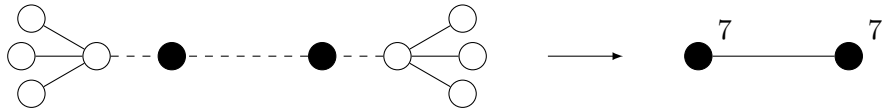

\begin{lemma}
	\label{lem:no22v}
	The graph $G$ cannot have adjacent 2-vertices.
\end{lemma}
\begin{proof}
	Suppose it did have vertices $u,v$ of degree 2 which are adjacent. 
	Let $X = \{u,v\}$, let $Y$ be the three edges incident to $u$ and $v$, and let $R$ be empty. (See Figure \ref{fig:no22v})
	\LemmaReduce
	\GsquaredChoose
	Therefore, $G$ cannot have adjacent 2-vertices.
\end{proof}

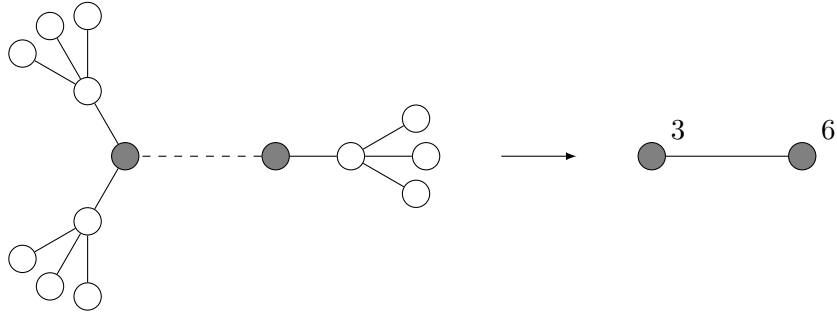
\begin{figure}[ht]
    \begin{center}
    \begin{tikzpicture}[node distance=2cm, every node/.style={draw=black, circle}]
        \node[fill=\recolored] (A) at (0*360/2:1) {};
        \node[fill=\recolored] (B) at (1*360/2:1) {};
        \node (Am) at ($(A) + (0:1)$) {};
        \node (Amm) at ($(Am) + (30:1)$) {};
        \node (Amn) at ($(Am) + (0:1)$) {};
        \node (Amo) at ($(Am) + (-30:1)$) {};
        \node (Bm) at ($(B) + (180+60:1)$) {};
        \node (Bmm) at ($(Bm) + (180+60+30:1)$) {};
        \node (Bmn) at ($(Bm) + (180+60+0:1)$) {};
        \node (Bmo) at ($(Bm) + (180+60-30:1)$) {};
        \node (Bn) at ($(B) + (180-60:1)$) {};
        \node (Bnm) at ($(Bn) + (180-60+30:1)$) {};
        \node (Bnn) at ($(Bn) + (180-60+0:1)$) {};
        \node (Bno) at ($(Bn) + (180-60-30:1)$) {};
        \node[fill=\recolored] (Az) at ($(A) + (GraphShift)$) {};
        \node[draw=none] (AzLabel) at ($(Az) + (NumShift)$) {6};
        \node[fill=\recolored] (Bz) at ($(B) + (GraphShift)$) {};
        \node[draw=none] (BzLabel) at ($(Bz) + (NumShift)$) {3};
        \draw[\edgeremoved] (A) -- (B);
        \draw (Am) -- (A);
        \draw (Amm) -- (Am) -- (Amn) (Am) -- (Amo);
        \draw (Bm) -- (B) -- (Bn);
        \draw (Bmm) -- (Bm) -- (Bmn) (Bm) -- (Bmo);
        \draw (Bnm) -- (Bn) -- (Bnn) (Bn) -- (Bno);
        \draw (Az) -- (Bz);
        \tikzArrow;
    \end{tikzpicture}
    \caption{A 2-vertex adjacent to a 3-vertex}
    \label{fig:no23v}
    \end{center}
\end{figure}

\begin{lemma}
	\label{lem:no23v}
	The graph $G$ cannot have a 2-vertex adjacent to a 3-vertex.
\end{lemma}
\begin{proof}
	Suppose it did have vertex $u$ of degree 2 and $v$ of degree 3 which are adjacent.
	Let $X$ be empty, let $Y = \{uv\}$, and let $R = \{u,v\}$. (See Figure \ref{fig:no23v})
	\LemmaReduce
	\GsquaredChoose
	Therefore, $G$ cannot have a 2-vertex adjacent to a 2-vertex.
\end{proof}

\begin{lemma}
	\label{lem:no33v}
	The graph $G$ cannot have a 3-vertex adjacent to a 3-vertex.
\end{lemma}
\begin{proof}
	Suppose it did have vertex $u$ of degree 3 and $v$ of degree 3 which are adjacent.
	Let $X$ be empty, let $Y = \{uv\}$, and let $R = \{u,v\}$.
	\LemmaReduce
	\GsquaredChoose
	Therefore, $G$ cannot have a 3-vertex adjacent to a 3-vertex.
\end{proof}

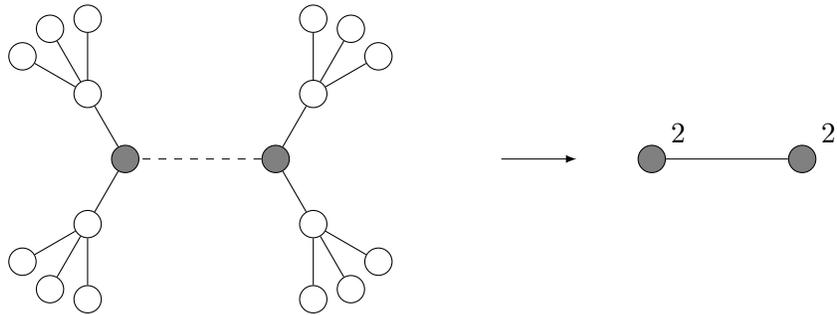
\begin{figure}[ht]
    \begin{center}
    \begin{tikzpicture}[node distance=2cm, every node/.style={draw=black, circle}]
        \node[fill=\recolored] (A) at (0*360/2:1) {};
        \node[fill=\recolored] (B) at (1*360/2:1) {};
        \node (Am) at ($(A) + (60:1)$) {};
        \node (Amm) at ($(Am) + (60+30:1)$) {};
        \node (Amn) at ($(Am) + (60+0:1)$) {};
        \node (Amo) at ($(Am) + (60-30:1)$) {};
        \node (An) at ($(A) + (-60:1)$) {};
        \node (Anm) at ($(An) + (-60+30:1)$) {};
        \node (Ann) at ($(An) + (-60+0:1)$) {};
        \node (Ano) at ($(An) + (-60-30:1)$) {};
        \node (Bm) at ($(B) + (180+60:1)$) {};
        \node (Bmm) at ($(Bm) + (180+60+30:1)$) {};
        \node (Bmn) at ($(Bm) + (180+60+0:1)$) {};
        \node (Bmo) at ($(Bm) + (180+60-30:1)$) {};
        \node (Bn) at ($(B) + (180-60:1)$) {};
        \node (Bnm) at ($(Bn) + (180-60+30:1)$) {};
        \node (Bnn) at ($(Bn) + (180-60+0:1)$) {};
        \node (Bno) at ($(Bn) + (180-60-30:1)$) {};
        \node[fill=\recolored] (Az) at ($(A) + (GraphShift)$) {};
        \node[draw=none] (AzLabel) at ($(Az) + (NumShift)$) {2};
        \node[fill=\recolored] (Bz) at ($(B) + (GraphShift)$) {};
        \node[draw=none] (BzLabel) at ($(Bz) + (NumShift)$) {2};
        \draw[\edgeremoved] (A) -- (B);
        \draw (Am) -- (A) -- (An);
        \draw (Amm) -- (Am) -- (Amn) (Am) -- (Amo);
        \draw (Anm) -- (An) -- (Ann) (An) -- (Ano);
        \draw (Bm) -- (B) -- (Bn);
        \draw (Bmm) -- (Bm) -- (Bmn) (Bm) -- (Bmo);
        \draw (Bnm) -- (Bn) -- (Bnn) (Bn) -- (Bno);
        \draw (Az) -- (Bz);
        \tikzArrow;
    \end{tikzpicture}
    \caption{Adjacent 3-vertices}
    \label{fig:no33v}
    \end{center}
\end{figure}

\begin{figure}[ht]
    \begin{center}
    \begin{tikzpicture}[node distance=2cm, every node/.style={draw=black, circle}]
        \node[fill=\vertremoved] (A) at (0*360/2:1) {};
        \node[fill=\recolored] (B) at (0,0) {};
        \node[fill=\vertremoved] (C) at (1*360/2:1) {};
        \node (Am) at ($(A) + (0:1)$) {};
        \node (Amm) at ($(Am) + (30:1)$) {};
        \node (Amn) at ($(Am) + (0:1)$) {};
        \node (Amo) at ($(Am) + (-30:1)$) {};
        \node (Bm) at ($(B) + (60:1)$) {};
        \node (Bn) at ($(B) + (120:1)$) {};
        \node (Bmm) at ($(Bm) + (60-30:1)$) {};
        \node (Bmn) at ($(Bm) + (60:1)$) {};
        \node (Bmo) at ($(Bm) + (60+30:1)$) {};
        \node (Bnm) at ($(Bn) + (1200-30:1)$) {};
        \node (Bnn) at ($(Bn) + (120:1)$) {};
        \node (Bno) at ($(Bn) + (120+30:1)$) {};
        \node (Cm) at ($(C) + (180:1)$) {};
        \node (Cmm) at ($(Cm) + (180+30:1)$) {};
        \node (Cmn) at ($(Cm) + (180+0:1)$) {};
        \node (Cmo) at ($(Cm) + (180-30:1)$) {};
        \node[fill=\vertremoved] (Az) at ($(A) + (GraphShift)$) {};
        \node[draw=none] (AzLabel) at ($(Az) + (NumShift)$) {6};
        \node[fill=\recolored] (Bz) at ($(B) + (GraphShift)$) {};
        \node[draw=none] (BzLabel) at ($(Bz) + (NumShift)$) {2};
        \node[fill=\vertremoved] (Cz) at ($(C) + (GraphShift)$) {};
        \node[draw=none] (CzLabel) at ($(Cz) + (NumShift)$) {6};
        \draw[\edgeremoved] (A) -- (B) -- (C);
        \draw[\edgeremoved] (Am) -- (A);
        \draw (Amm) -- (Am) -- (Amn) (Am) -- (Amo);
        \draw (Bm) -- (B) -- (Bn);
        \draw (Bmm) -- (Bm) -- (Bmn) (Bm) -- (Bmo);
        \draw (Bnm) -- (Bn) -- (Bnn) (Bn) -- (Bno);
        \draw[dashed] (Cm) -- (C);
        \draw (Cmm) -- (Cm) -- (Cmn) (Cm) -- (Cmo);
        \draw (Az) -- (Bz) -- (Cz);
        \draw (Az) edge[bend left] (Cz);
        \tikzArrow;
    \end{tikzpicture}
    \caption{2-vertices at distance at most two}
    \label{fig:no242v}
    \end{center}
\end{figure}
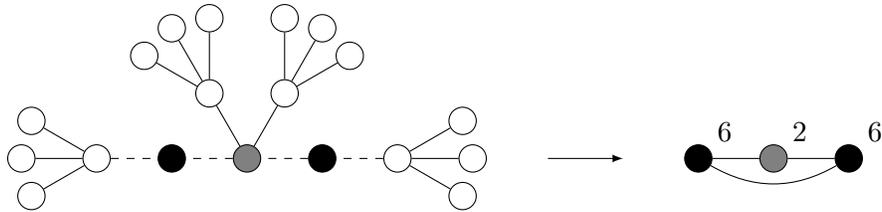

\begin{lemma}
	\label{lem:no242v}
	The graph $G$ cannot have a 2-vertex distance at most two away from another 2-vertex.
\end{lemma}
\begin{proof}
	Suppose it did have vertices $u,v,w$ such that $u,w$ are 2-vertices, and $uvw$ is a path. Note that by Lemma \ref{lem:no23v}, $v$ must be a 4-vertex.
	Let $X = \{u,w\}$, let $Y$ be the four edges incident to $u,w$, and let $R = \{v\}$. (See Figure \ref{fig:no242v})
	\LemmaReduce
	\GsquaredChoose
	Therefore, $G$ cannot have a 2-vertex distance two away from another 2-vertex.
\end{proof}

\begin{figure}[ht]
    \begin{center}
    \begin{tikzpicture}[node distance=2cm, every node/.style={draw=black, circle}]
        \node[fill=\vertremoved] (A) at (0*360/2:1) {};
        \node[fill=\recolored] (B) at (0,0) {};
        \node[fill=\recolored] (C) at (1*360/2:1) {};
        \node (Am) at ($(A) + (0:1)$) {};
        \node (Amm) at ($(Am) + (30:1)$) {};
        \node (Amn) at ($(Am) + (0:1)$) {};
        \node (Amo) at ($(Am) + (-30:1)$) {};
        \node (Bm) at ($(B) + (60:1)$) {};
        \node (Bn) at ($(B) + (120:1)$) {};
        \node (Bmm) at ($(Bm) + (60-30:1)$) {};
        \node (Bmn) at ($(Bm) + (60:1)$) {};
        \node (Bmo) at ($(Bm) + (60+30:1)$) {};
        \node (Bnm) at ($(Bn) + (1200-30:1)$) {};
        \node (Bnn) at ($(Bn) + (120:1)$) {};
        \node (Bno) at ($(Bn) + (120+30:1)$) {};
        \node (Cm) at ($(C) + (180-30:1)$) {};
        \node (Cmm) at ($(Cm) + (180-30+30:1)$) {};
        \node (Cmn) at ($(Cm) + (180-30+0:1)$) {};
        \node (Cmo) at ($(Cm) + (180-30-30:1)$) {};
        \node (Cn) at ($(C) + (180+30:1)$) {};
        \node (Cnm) at ($(Cn) + (180+30+30:1)$) {};
        \node (Cnn) at ($(Cn) + (180+30+0:1)$) {};
        \node (Cno) at ($(Cn) + (180+30-30:1)$) {};
        \node[fill=\vertremoved] (Az) at ($(A) + (GraphShift)$) {};
        \node[draw=none] (AzLabel) at ($(Az) + (NumShift)$) {6};
        \node[fill=\recolored] (Bz) at ($(B) + (GraphShift)$) {};
        \node[draw=none] (BzLabel) at ($(Bz) + (NumShift)$) {1};
        \node[fill=\recolored] (Cz) at ($(C) + (GraphShift)$) {};
        \node[draw=none] (CzLabel) at ($(Cz) + (NumShift)$) {2};
        \draw[dashed] (A) -- (B);
        \draw (B) -- (C);
        \draw[\edgeremoved] (Am) -- (A);
        \draw (Amm) -- (Am) -- (Amn) (Am) -- (Amo);
        \draw (Bm) -- (B) -- (Bn);
        \draw (Bmm) -- (Bm) -- (Bmn) (Bm) -- (Bmo);
        \draw (Bnm) -- (Bn) -- (Bnn) (Bn) -- (Bno);
        \draw (Cm) -- (C) -- (Cn);
        \draw (Cmm) -- (Cm) -- (Cmn) (Cm) -- (Cmo);
        \draw (Cnm) -- (Cn) -- (Cnn) (Cn) -- (Cno);
        \draw (Az) -- (Bz) -- (Cz);
        \draw (Az) edge[bend left] (Cz);
        \tikzArrow;
    \end{tikzpicture}
    \caption{A 2-vertex at distance at most two from a 3-vertex}
    \label{fig:no243v}
    \end{center}
\end{figure}
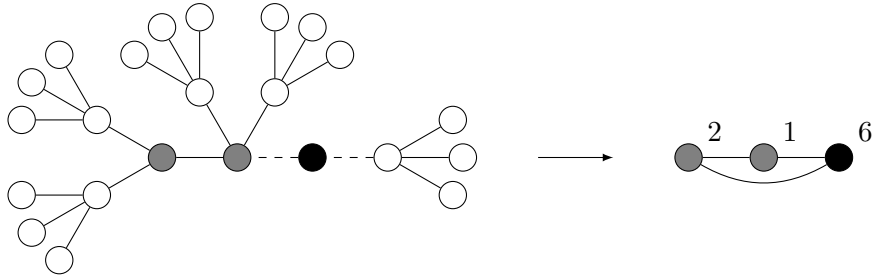

\begin{lemma}
	\label{lem:no243v}
	The graph $G$ cannot have a 2-vertex distance at most two away from a 3-vertex.
\end{lemma}
\begin{proof}
	Suppose it did have vertices $u,v,w$ such that $u$ is a 2-vertex, $w$ is a 3-vertex, and $uvw$ is a path. Note that by Lemma \ref{lem:no23v}, $v$ must be a 4-vertex.
	Let $X = \{u\}$, let $Y$ be the two edges incident to $u$, and let $R = \{v,w\}$. (See Figure \ref{fig:no243v})
	\LemmaReduce
	\GsquaredChoose
	Therefore, $G$ cannot have a 2-vertex distance two away from a 3-vertex.
\end{proof}

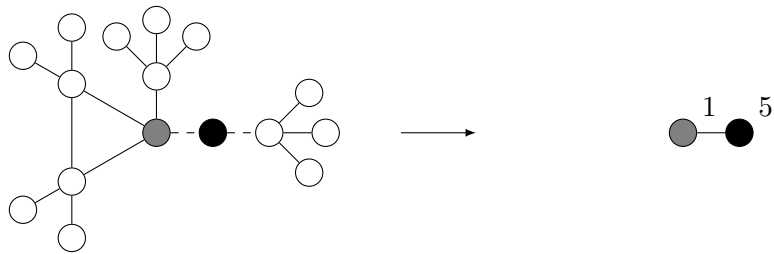
\begin{figure}[ht]
    \begin{center}
    \begin{tikzpicture}[node distance=2cm, every node/.style={draw=black, circle}]
        \node[fill=\recolored] (A) at (0*360/3:0.75) {};
        \node (B) at (1*360/3:0.75) {};
        \node (C) at (2*360/3:0.75) {};
        \node[fill=\vertremoved] (D) at (0:1.5) {};
        \node (Am) at ($(A) + (90:0.75)$) {};
        \node (Amm) at ($(Am) + (90+45:0.75)$) {};
        \node (Amn) at ($(Am) + (90+0:0.75)$) {};
        \node (Amo) at ($(Am) + (90-45:0.75)$) {};
        \node (Bm) at ($(B) + (1*360/3-30:0.75)$) {};
        \node (Bn) at ($(B) + (1*360/3+30:0.75)$) {};
        \node (Cm) at ($(C) + (2*360/3-30:0.75)$) {};
        \node (Cn) at ($(C) + (2*360/3+30:0.75)$) {};
        \node (Dm) at ($(D) + (0:0.75)$) {};
        \node (Dmm) at ($(Dm) + (0+45:0.75)$) {};
        \node (Dmn) at ($(Dm) + (0+0:0.75)$) {};
        \node (Dmo) at ($(Dm) + (0-45:0.75)$) {};
        \node[fill=\recolored] (Az) at ($(A) + (GraphShift)$) {};
        \node[draw=none] (AzLabel) at ($(Az) + (NumShift)$) {1};
        \node[fill=\vertremoved] (Dz) at ($(D) + (GraphShift)$) {};
        \node[draw=none] (DzLabel) at ($(Dz) + (NumShift)$) {5};
        \draw (A) -- (B) -- (C) -- (A);
        \draw[\edgeremoved] (A) -- (D);
        \draw (Am) -- (A);
        \draw (Amm) -- (Am) -- (Amn) (Am) -- (Amo);
        \draw (Bm) -- (B) -- (Bn);
        \draw (Cm) -- (C) -- (Cn);
        \draw[\edgeremoved] (D) -- (Dm);
        \draw (Dmm) -- (Dm) -- (Dmn) (Dm) -- (Dmo);
        \draw (Az) -- (Dz);
        \tikzArrow;
    \end{tikzpicture}
    \caption{A 2-vertex adjacent to a vertex of a 3-face}
    \label{fig:no2v-3f}
    \end{center}
\end{figure}

\begin{lemma}
    \label{lem:no2v-3f}
    The graph $G$ cannot have a 2-vertex adjacent to one of the vertices of a 3-face.
\end{lemma}
\begin{proof}
    Suppose it did have vertices $u,v$ such that $u$ is a 2-vertex, $v$ is incident to a 3-face $f$ and $u,v$ are adjacent. Note that by Lemmas \ref{lem:no22v}, \ref{lem:no23v}, \ref{lem:no242v}, and \ref{lem:no243v}, $v$ and all the vertices incident to the 3-face must be 4-vertices.
    Let $X = \{u\}$, let $Y$ be the two edges incident to $u$, and let $R = \{v\}$. (See Figure \ref{fig:no2v-3f})
    \LemmaReduce
    \GsquaredChoose
    Therefore, $G$ cannot have a 2-vertex adjacent to one of the vertices of a 3-face.
\end{proof}

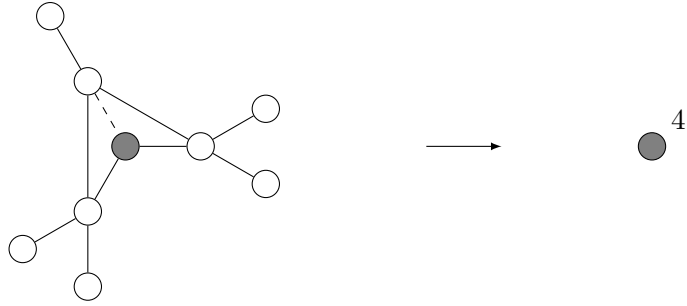
\begin{figure}[ht]
    \begin{center}
    \begin{tikzpicture}[node distance=2cm, every node/.style={draw=black, circle}]
        \node[fill=\recolored] (A) at (0,0) {};
        \node (B) at (0*360/3:1) {};
        \node (C) at (1*360/3:1) {};
        \node (D) at (2*360/3:1) {};
        \node (Bm) at ($(B) + (0*360/3-30:1)$) {};
        \node (Bn) at ($(B) + (0*360/3+30:1)$) {};
        \node (Cm) at ($(C) + (1*360/3-0:1)$) {};
        \node (Dm) at ($(D) + (2*360/3-30:1)$) {};
        \node (Dn) at ($(D) + (2*360/3+30:1)$) {};
        \node[fill=\recolored] (Az) at ($(A) + (GraphShift)$) {};
        \node[draw=none] (AzLabel) at ($(Az) + (NumShift)$) {4};
        \draw (B) -- (C) -- (D);
        \draw (B) -- (A) -- (D);
        \draw[\edgeremoved] (C) -- (A);
        \draw (Bm) -- (B) -- (Bn);
        \draw (Cm) -- (C);
        \draw (Dm) -- (D) -- (Dn);
        \tikzArrow;
    \end{tikzpicture}
    \caption{A 3-vertex incident to two 3-faces}
    \label{fig:no3v_33f}
    \end{center}
\end{figure}

\begin{lemma}
    \label{lem:no3v_33f}
    The graph $G$ can have no 3-vertex incident to two 3-faces.
\end{lemma}
\begin{proof}
    Suppose it did have a 3-vertex $v$ which is incident to two 3-faces $f,g$, then those 3-faces must share an edge $vu$. 
    Let $X$ be empty, $Y = \{ vu \}$ and $R = \{v\}$. (See Figure \ref{fig:no3v_33f})
    \LemmaReduce
    \GsquaredChoose
    Therefore, $G$ cannot have a 3-vertex incident to two 3-faces.
\end{proof}

\begin{lemma}
	\label{lem:no333f}
	The graph $G$ can have no 3-face sharing two edges with 3-faces.
\end{lemma}
\begin{proof}
	If the graph $G$ contains a 3-face $f$ which shares two edges $uv,vw$ with the same 3-face $g$, then the vertex $v$ must be a 2-vertex incident to a 3-face, which by Lemma \ref{lem:no2v3f} cannot happen.

	If the graph $G$ contains a 3-face which shares two edges with two distinct 3-faces which do not share an edge themselves, then the outside cycle around the three faces is a 5-cycle, which is forbidden. Therefore, this cannot happen.

	Suppose the graph $G$ contains a 3-face which shares two edges $uv,vw$ with two distinct 3-faces which also share an edge $vx$. Then $v$ is a 3-vertex which is incident to three 3-faces, which by Lemma \ref{lem:no3v_33f} cannot happen.

	Since none of these possibilities can happen, $G$ cannot have a 3-face sharing two edges with 3-faces.
\end{proof}

\begin{lemma}
	\label{lem:no34f}
	The graph $G$ can have no 3-face sharing an edge with a 4-face.
\end{lemma}
\begin{proof}
	If the graph $G$ contains a 3-face which shares two edges $uv, vw$ with the same 4-face, then the vertex $v$ must be a 2-vertex incident to a 3-face, which by Lemma \ref{lem:no2v3f} cannot happen.

	If the graph $G$ contains a 3-face which shares a single edge with a 4-face, then the outside cycle around the two faces is a 5-cycle, which is forbidden. Therefore, this cannot happen.
\end{proof}

\begin{lemma}
	\label{lem:no3v_44f}
	The graph $G$ can have no 3-vertex incident to two 4-faces.
\end{lemma}
\begin{proof}
	If the graph $G$ contains a 3-vertex $v$ incident to two 4-faces, then those two 4-faces must share an edge $vu$. 
	Let $X$ be empty, $Y = \{vu\}$, and let $R = \{v\}$.
	\LemmaReduce
	\GsquaredChoose
	Therefore, $G$ cannot have a 3-vertex incident to two 4-faces.
\end{proof}

\begin{figure}[ht]
    \begin{center}
    \begin{tikzpicture}[node distance=2cm, every node/.style={draw=black, circle}]
        \node[fill=\recolored] (A) at (0,0) {};
        \node (B) at (0*360/4:1) {};
        \node (C) at (1*360/4:1) {};
        \node (D) at (2*360/4:1) {};
        \node (BC) at (0.5*360/4:\roottwo) {};
        \node (CD) at (1.5*360/4:\roottwo) {};
        \node (Bm) at ($(B) + (0*360/4-30:1)$) {};
        \node (Bn) at ($(B) + (0*360/4+30:1)$) {};
        \node (Cm) at ($(C) + (1*360/4-0:1)$) {};
        \node (Dm) at ($(D) + (2*360/4-30:1)$) {};
        \node (Dn) at ($(D) + (2*360/4+30:1)$) {};
        \node[fill=\recolored] (Az) at ($(A) + (GraphShift)$) {};
        \node[draw=none] (AzLabel) at ($(Az) + (NumShift)$) {2};
        \draw (B) -- (BC) -- (C) -- (CD) -- (D);
        \draw (B) -- (A) -- (D);
        \draw[\edgeremoved] (C) -- (A);
        \draw (Bm) -- (B) -- (Bn);
        \draw (Cm) -- (C);
        \draw (Dm) -- (D) -- (Dn);
        \tikzArrow;
    \end{tikzpicture}
    \caption{A 3-vertex incident to two 4-faces}
    \label{fig:no3v_44f}
    \end{center}
\end{figure}
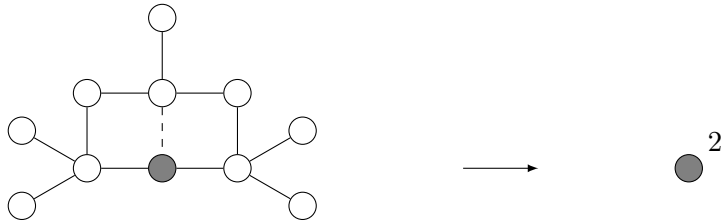

\begin{figure}[ht]
    \begin{center}
    \begin{tikzpicture}[node distance=2cm, every node/.style={draw=black, circle}]
        \node (A) at (0*360/3:1) {};
        \node[fill=\recolored] (B) at (1*360/3:1) {};
        \node[fill=\recolored] (C) at (2*360/3:1) {};
        \node (D) at (180:2) {};
        \node (Am) at ($(A) + (0*360/3:1)$) {};
        \node (Bm) at ($(B) + (1*360/3:1)$) {};
        \node (Bmm) at ($(Bm) + (1*360/3-30:1)$) {};
        \node (Bmn) at ($(Bm) + (1*360/3+0:1)$) {};
        \node (Bmo) at ($(Bm) + (1*360/3+30:1)$) {};
        \node (Cm) at ($(C) + (2*360/3:1)$) {};
        \node (Cmm) at ($(Cm) + (2*360/3-30:1)$) {};
        \node (Cmn) at ($(Cm) + (2*360/3+0:1)$) {};
        \node (Cmo) at ($(Cm) + (2*360/3+30:1)$) {};
        \node (Dm) at ($(D) + (180-30:1)$) {};
        \node (Dn) at ($(D) + (180+30:1)$) {};
        \node[fill=\recolored] (Bz) at ($(B) + (GraphShift)$) {};
        \node[draw=none] (BzLabel) at ($(Bz) + (NumShift)$) {2};
        \node[fill=\recolored] (Cz) at ($(C) + (GraphShift)$) {};
        \node[draw=none] (CzLabel) at ($(Cz) + (NumShift)$) {2};
        \draw[\edgeremoved] (B) -- (C);
        \draw (C) -- (A) -- (B);
        \draw (C) -- (D) -- (B);
        \draw (Am) -- (A);
        \draw (Bm) -- (B);
        \draw (Bmm) -- (Bm) -- (Bmn) (Bm) -- (Bmo);
        \draw (Cm) -- (C);
        \draw (Cmm) -- (Cm) -- (Cmn) (Cm) -- (Cmo);
        \draw (Dm) -- (D) -- (Dn);
        \draw (Bz) -- (Cz);
        \tikzArrow;
    \end{tikzpicture}
    \caption{A 3-vertex incident a 3-face sharing an edge with another 3-face}
    \label{fig:no3v3f3f}
    \end{center}
\end{figure}
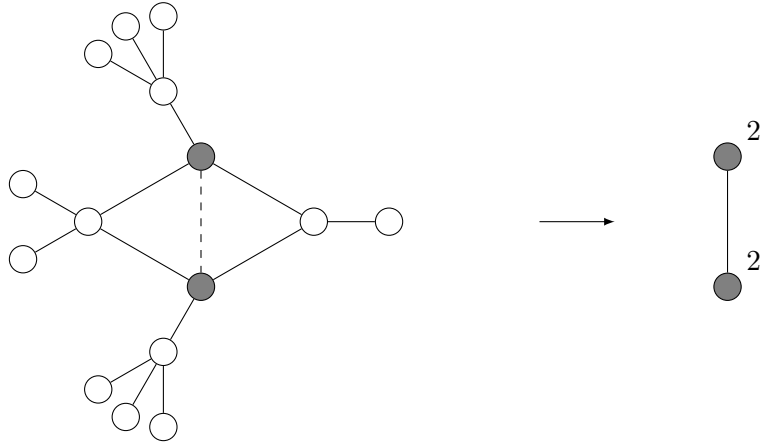

\begin{lemma}
    \label{lem:no3v3f3f}
    The graph $G$ can have no 3-vertex incident to a 3-face which shares an edge with another 3-face.
\end{lemma}
\begin{proof}
    Suppose it did have a 3-vertex $v$ incident to a 3-face $f$ which shares an edge with another 3-face $g$. By Lemma \ref{lem:no3v_33f}, $v$ cannot be incident to the 3-face $g$ also, therefore, the shared edge between the faces $f,g$ must not be incident to $v$. Called the shared edge $uw$. Therefore, $f$ is incident to the vertices $v,u,w$, and $g$ is incident to $u,w$ and a third vertex. 
    Let $X$ be empty, $Y = \{uw\}$ and let $R = \{u,w\}$. (See Figure \ref{fig:no3v3f3f})
    \LemmaReduce
    \GsquaredChoose
    Therefore, $G$ can have no 3-vertex incident to a 3-face which shares an edge with another 3-face.
\end{proof}

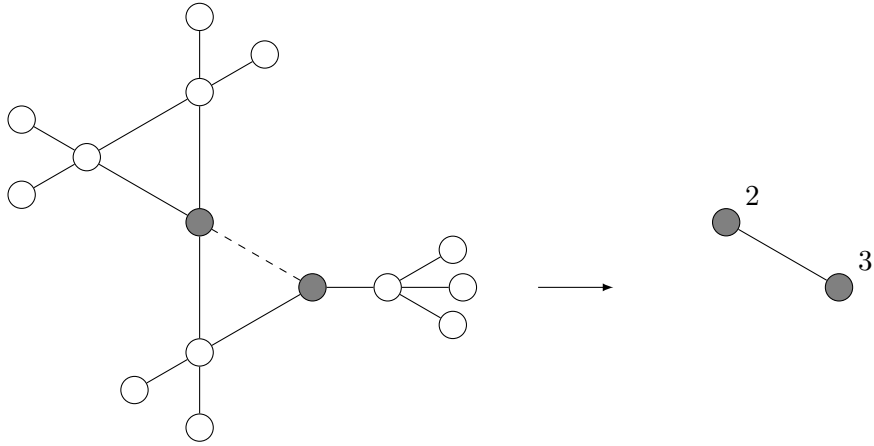
\begin{figure}[ht]
    \begin{center}
    \begin{tikzpicture}[node distance=2cm, every node/.style={draw=black, circle}]
        \node[fill=\recolored] (A) at (0*360/3:1) {};
        \node[fill=\recolored] (B) at (1*360/3:1) {};
        \node (C) at (2*360/3:1) {};
        \node (D) at ($2*(B) - (A)$) {};
        \node (E) at ($2*(B) - (C)$) {};
        \node (Am) at ($(A) + (0*360/3:1)$) {};
        \node (Amm) at ($(Am) + (0*360/3-30:1)$) {};
        \node (Amn) at ($(Am) + (0*360/3+0:1)$) {};
        \node (Amo) at ($(Am) + (0*360/3+30:1)$) {};
        \node (Cm) at ($(C) + (2*360/3-30:1)$) {};
        \node (Cn) at ($(C) + (2*360/3+30:1)$) {};
        \node (Dm) at ($(D) + (180-30:1)$) {};
        \node (Dn) at ($(D) + (180+30:1)$) {};
        \node (Em) at ($(E) + (1*360/6-30:1)$) {};
        \node (En) at ($(E) + (1*360/6+30:1)$) {};
        \node[fill=\recolored] (Az) at ($(A) + (GraphShift)$) {};
        \node[draw=none] (AzLabel) at ($(Az) + (NumShift)$) {3};
        \node[fill=\recolored] (Bz) at ($(B) + (GraphShift)$) {};
        \node[draw=none] (BzLabel) at ($(Bz) + (NumShift)$) {2};
        \draw[\edgeremoved] (A) -- (B);
        \draw (B) -- (C) -- (A);
        \draw (D) -- (E) -- (B) -- (D);
        \draw (Am) -- (A);
        \draw (Amm) -- (Am) -- (Amn) (Am) -- (Amo);
        \draw (Cm) -- (C) -- (Cn);
        \draw (Dm) -- (D) -- (Dn);
        \draw (Em) -- (E) -- (En);
        \draw (Az) -- (Bz);
        \tikzArrow;
    \end{tikzpicture}
    \caption{A 3-vertex incident a 3-face sharing a vertex with another 3-face}
    \label{fig:no3v3f-3f}
    \end{center}
\end{figure}

\begin{lemma}
    \label{lem:no3v3f-3f}
    The graph $G$ can have no 3-vertex incident to a 3-face which shares a vertex with another 3-face.
\end{lemma}
\begin{proof}
    Suppose it does have a 3-vertex $v$ incident to a 3-face $f$ which shares a vertex with another 3-face $g$. By Lemma \ref{lem:no3v_33f}, $v$ cannot also be incident to the 3-face $g$. Therefore, the vertex shared between faces $f$ and $g$ must be some other vertex we will call $u$.
    Let $X$ be empty, let $Y = \{uv\}$, and let $R = \{u,v\}$. (See Figure \ref{fig:no3v3f-3f})
    \LemmaReduce
    \GsquaredChoose
    Therefore, $G$ can have no 3-vertex incident to a 3-face which shares a vertex with another 3-face.
\end{proof}

\begin{figure}[ht]
    \begin{center}
    \begin{tikzpicture}[node distance=2cm, every node/.style={draw=black, circle}]
        \node[fill=\recolored] (A) at (0*360/3:1) {};
        \node[fill=\recolored] (B) at (1*360/3:1) {};
        \node (C) at (2*360/3:1) {};
        \node (D) at (0:2) {};
        \node (Am) at ($(A) + (0*360/3-80:1)$) {};
        \node (Amm) at ($(Am) + (0*360/3-80-30:1)$) {};
        \node (Amn) at ($(Am) + (0*360/3-80+0:1)$) {};
        \node (Amo) at ($(Am) + (0*360/3-80+30:1)$) {};
        \node (Bm) at ($(B) + (1*360/3:1)$) {};
        \node (Bmm) at ($(Bm) + (1*360/3-30:1)$) {};
        \node (Bmn) at ($(Bm) + (1*360/3+0:1)$) {};
        \node (Bmo) at ($(Bm) + (1*360/3+30:1)$) {};
        \node (Cm) at ($(C) + (2*360/3-30:1)$) {};
        \node (Cn) at ($(C) + (2*360/3+30:1)$) {};
        \node (Dm) at ($(D) + (0-30:1)$) {};
        \node (Dn) at ($(D) + (0+30:1)$) {};
        \node[fill=\recolored] (Az) at ($(A) + (GraphShift)$) {};
        \node[draw=none] (AzLabel) at ($(Az) + (NumShift)$) {1};
        \node[fill=\recolored] (Bz) at ($(B) + (GraphShift)$) {};
        \node[draw=none] (BzLabel) at ($(Bz) + (NumShift)$) {3};
        \draw[\edgeremoved] (A) -- (B);
        \draw (B) -- (C) -- (A);
        \draw (A) -- (D);
        \draw (Am) -- (A);
        \draw (Amm) -- (Am) -- (Amn) (Am) -- (Amo);
        \draw (Bm) -- (B);
        \draw (Bmm) -- (Bm) -- (Bmn) (Bm) -- (Bmo);
        \draw (Cm) -- (C) -- (Cn);
        \draw (Dm) -- (D) -- (Dn);
        \draw (Az) -- (Bz);
        \tikzArrow;
    \end{tikzpicture}
    \caption{A 3-vertex adjacent to a vertex incident to a 3-face incident to another 3-vertex}
    \label{fig:no3v-3f3v}
    \end{center}
\end{figure}
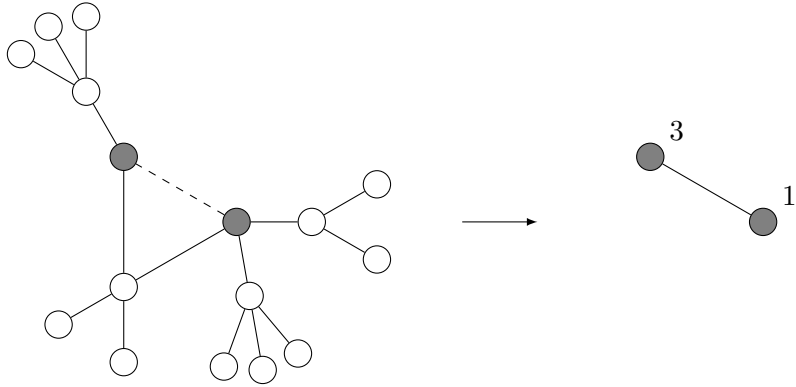

\begin{lemma}
    \label{lem:no3v-3f3v}
    The graph $G$ can have no 3-vertex adjacent to a vertex incident to a 3-face which is incident to another 3-vertex.
\end{lemma}
\begin{proof}
    Suppose it does have a 3-vertex $v$ which is adjacent to a vertex $u$ incident to a 3-face $f$ that is incident to another 3-vertex. By Lemma \ref{lem:no33v}, $u$ cannot be a 3-vertex, therefore, the 3-vertex incident to $f$ must be some other vertex $w$, and the vertex $u$ must be a 4-vertex. 
    Let $X$ be empty, $Y = \{uw\}$, and let $R = \{u,w\}$. (See Figure \ref{fig:no3v-3f3v})
    \LemmaReduce
    \GsquaredChoose
    Therefore, $G$ can have no 3-vertex adjacent to a vertex incident to a 3-face which is incident to another 3-vertex.
\end{proof}

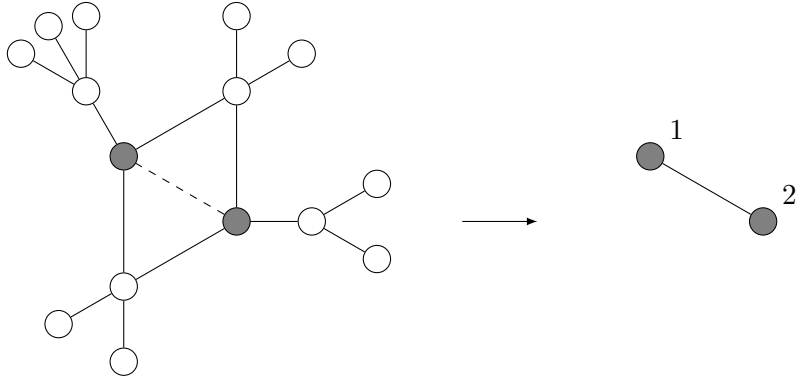
\begin{figure}[ht]
    \begin{center}
    \begin{tikzpicture}[node distance=2cm, every node/.style={draw=black, circle}]
        \node[fill=\recolored] (A) at (0*360/3:1) {};
        \node[fill=\recolored] (B) at (1*360/3:1) {};
        \node (C) at (2*360/3:1) {};
        \node (D) at (360/6:2) {};
        \node (Am) at ($(A) + (0*360/3:1)$) {};
        \node (Amm) at ($(Am) + (0*360/3-30:1)$) {};
        \node (Amn) at ($(Am) + (0*360/3+30:1)$) {};
        \node (Bm) at ($(B) + (1*360/3:1)$) {};
        \node (Bmm) at ($(Bm) + (1*360/3-30:1)$) {};
        \node (Bmn) at ($(Bm) + (1*360/3+0:1)$) {};
        \node (Bmo) at ($(Bm) + (1*360/3+30:1)$) {};
        \node (Cm) at ($(C) + (2*360/3-30:1)$) {};
        \node (Cn) at ($(C) + (2*360/3+30:1)$) {};
        \node (Dm) at ($(D) + (360/6-30:1)$) {};
        \node (Dn) at ($(D) + (360/6+30:1)$) {};
        \node[fill=\recolored] (Az) at ($(A) + (GraphShift)$) {};
        \node[draw=none] (AzLabel) at ($(Az) + (NumShift)$) {2};
        \node[fill=\recolored] (Bz) at ($(B) + (GraphShift)$) {};
        \node[draw=none] (BzLabel) at ($(Bz) + (NumShift)$) {1};
        \draw[\edgeremoved] (A) -- (B);
        \draw (B) -- (C) -- (A);
        \draw (B) -- (D) -- (A);
        \draw (Am) -- (A);
        \draw (Amm) -- (Am) -- (Amn);
        \draw (Bm) -- (B);
        \draw (Bmm) -- (Bm) -- (Bmn) (Bm) -- (Bmo);
        \draw (Cm) -- (C) -- (Cn);
        \draw (Dm) -- (D) -- (Dn);
        \draw (Az) -- (Bz);
        \tikzArrow;
    \end{tikzpicture}
    \caption{A 3-vertex adjacent to one of the endpoints of the shared edge between two 3-faces}
    \label{fig:no3v-m3f3f}
    \end{center}
\end{figure}

\begin{lemma}
    \label{lem:no3v-m3f3f}
    The graph $G$ can have no $3^-$-vertex adjacent to one of the endpoints of the shared edge between two 3-faces.
\end{lemma}
\begin{proof}
    Suppose it does have a $3^-$-vertex $v$ which is adjacent to a vertex $u$ such that the edge $uw$ is the shared edge between two 3-faces $f,g$.
    We know that $w \ne v$ by Lemma \ref{lem:no3v_33f} if $v$ is a 3-vertex, or by Lemma \ref{lem:no2v3f} if $v$ is a 2-vertex. Without loss of generality, we can now assume that $v$ is a 3-vertex.
    Let $X$ be empty, let $Y = \{uw\}$, and let $R = \{u,w\}$. (See Figure \ref{fig:no3v-m3f3f})
    \LemmaReduce
    \GsquaredChoose
    Therefore, $G$ can have no 3-vertex adjacent to one of the endpoints of the shared edge between two 3-faces.
\end{proof}

\begin{figure}[ht]
    \begin{center}
    \begin{tikzpicture}[node distance=2cm, every node/.style={draw=black, circle}]
        \node[fill=\recolored] (A) at (0*360/3:1) {};
        \node[fill=\recolored] (B) at (1*360/3:1) {};
        \node (C) at (2*360/3:1) {};
        \node (D) at (360/6:2) {};
        \node[fill=\recolored] (E) at ($(A) + (0:1)$) {};
        \node[fill=\recolored] (F) at ($(E) + (0:1)$) {};
        \node (Bm) at ($(B) + (1*360/3:1)$) {};
        \node (Bmm) at ($(Bm) + (1*360/3-30:1)$) {};
        \node (Bmn) at ($(Bm) + (1*360/3+0:1)$) {};
        \node (Bmo) at ($(Bm) + (1*360/3+30:1)$) {};
        \node (Cm) at ($(C) + (2*360/3-30:1)$) {};
        \node (Cn) at ($(C) + (2*360/3+30:1)$) {};
        \node (Dm) at ($(D) + (360/6-30:1)$) {};
        \node (Dn) at ($(D) + (360/6+30:1)$) {};
        \node (Em) at ($(E) + (-90-30:1)$) {};
        \node (Emm) at ($(Em) + (-90-30-30:1)$) {};
        \node (Emn) at ($(Em) + (-90-30+0:1)$) {};
        \node (Emo) at ($(Em) + (-90-30+30:1)$) {};
        \node (En) at ($(E) + (-90+30:1)$) {};
        \node (Enm) at ($(En) + (-90+30-30:1)$) {};
        \node (Enn) at ($(En) + (-90+30+0:1)$) {};
        \node (Eno) at ($(En) + (-90+30+30:1)$) {};
        \node (Fm) at ($(F) + (90:1)$) {};
        \node (Fmm) at ($(Fm) + (90-30:1)$) {};
        \node (Fmn) at ($(Fm) + (90+0:1)$) {};
        \node (Fmo) at ($(Fm) + (90+30:1)$) {};
        \node[fill=\recolored] (Az) at ($(A) + (GraphShift)$) {};
        \node[draw=none] (AzLabel) at ($(Az) + (NumShift)$) {2};
        \node[fill=\recolored] (Bz) at ($(B) + (GraphShift)$) {};
        \node[draw=none] (BzLabel) at ($(Bz) + (NumShift)$) {3};
        \node[fill=\recolored] (Ez) at ($(E) + (GraphShift)$) {};
        \node[draw=none] (EzLabel) at ($(Ez) + (NumShift)$) {1};
        \node[fill=\recolored] (Fz) at ($(F) + (GraphShift)$) {};
        \node[draw=none] (FzLabel) at ($(Fz) + (NumShift)$) {6};
        \draw[\edgeremoved] (A) -- (B);
        \draw (B) -- (C) -- (A);
        \draw (B) -- (D) -- (A);
        \draw (A) -- (E) -- (F);
        \draw (Bm) -- (B);
        \draw (Bmm) -- (Bm) -- (Bmn) (Bm) -- (Bmo);
        \draw (Cm) -- (C) -- (Cn);
        \draw (Dm) -- (D) -- (Dn);
        \draw (Em) -- (E) -- (En);
        \draw (Emm) -- (Em) -- (Emn) (Em) -- (Emo);
        \draw (Enm) -- (En) -- (Enn) (En) -- (Eno);
        \draw (Em) -- (E) -- (En);
        \draw (Emm) -- (Em) -- (Emn) (Em) -- (Emo);
        \draw (Fm) -- (F);
        \draw (Fmm) -- (Fm) -- (Fmn) (Fm) -- (Fmo);
        \draw (Bz) -- (Az) -- (Ez) -- (Fz);
        \draw (Bz) edge[bend left] (Ez);
        \draw (Az) edge[bend right] (Fz);
        \tikzArrow;
    \end{tikzpicture}
    \caption{A 2-vertex distance at most two from one of the endpoints of the shared edge between two 3-faces}
    \label{fig:no2v--m3f3f}
    \end{center}
\end{figure}
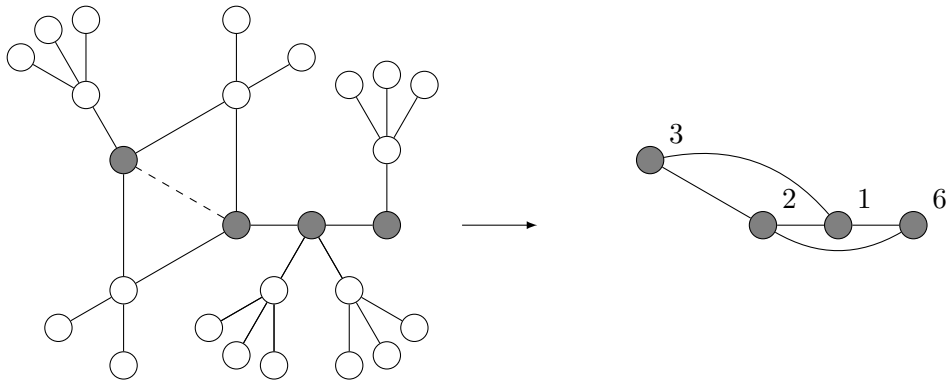

\begin{lemma}
    \label{lem:no2v--m3f3f}
    The graph $G$ can have no 2-vertex distance at most two away from the midpoint of one of the endpoints of the shared edge between two 3-faces.
\end{lemma}
\begin{proof}
    Suppose it does have a path $vuwx$ such that $v$ is of degree 2, and $wx$ is the shared edge between two 3-faces. By Lemma \ref{lem:no3v-m3f3f}, $u$ must be of degree 4.
    Let $X$ be empty, let $Y = \{wx\}$, and let $R = \{v,u,w,x\}$. (See Figure \ref{fig:no2v--m3f3f})
    \LemmaReduce
    The subgraph of $G^2$ induced by $X \cup R$ in this general instance is not complete, it is missing one edge. However, even if that edge were present and the $f$-values remained the same, the resulting graph would still be $f$-choosable, and so the conditions of Lemma \ref{lem:reduceWegner4} are still satisfied.
    \GsquaredChoose
    Therefore, $G$ can have no 2-vertex distance at most two away from the midpoint of one of the endpoints of the shared edge between two 3-faces
\end{proof}

\section{Discharging}
\label{sec:color_discharge}

In this proof, we will use the technique of \emph{intermediate discharging}. This is a technique where the main discharging is completed by a series of smaller steps within the proof. Usually this is done in order to clarify the proof.

\begin{proof}{(Theorem \ref{thm:wegner4no5})}
	To each vertex of $G$, give a charge equal to its degree minus four, and to each face of $G$ give a charge equal to its length minus four. This is commonly known as balanced charging. Let $c$ be this charging function. By using the handshaking lemma and Euler's formula (which we can use by Lemma \ref{lem:conn}), we can derive the total charge on $G$:
	\begin{align*}
		\sum_{v \in V(G)} c(v) + \sum_{f \in f(G)} c(f) &= \sum_{v \in V(G)} (d(v) - 4) + \sum_{f \in F(G)} (l(f) - 4) \\
		&= \left(\sum_{v \in V(G)} d(v)\right) - 4|V(G)| + \left(\sum_{f \in F(G)} l(f)\right) - 4|F(G)| \\
		&= 2|E(G)| - 4|V(G)| + 2|E(G)| - 4|F(G)| \\
		&= -4(|V(G)| - |E(G)| + |F(G)|) \\
		&= -4 \cdot 2 \\
		&= -8
	\end{align*}
	Thus, the total initial charge is negative.

	We then redistribute the charge by the following discharging rules:
	\begin{enumerate}[R1:]
		\item A 2-vertex takes 1 charge from each incident $6^+$-face. \label{rule:2v}
		\item A 3-vertex takes $\frac{1}{2}$ charge from each incident $6^+$-face. \label{rule:3v}
		\item A 3-face incident to at least one 3-face takes $\frac{1}{2}$ charge from each incident $6^+$-face. \label{rule:3f3f}
		\item A 3-face not incident to any 3-faces takes $\frac{1}{3}$ charge from each incident $6^+$-face. \label{rule:3f}
	\end{enumerate}
	In the above rules, if a vertex is incident to a face in multiple ways, and it would take charge from that face according to the above rules, we want it to take charge from that face multiple times. Similarly, we want the same to apply for faces. For example, a 2-vertex $v$ is generally incident to two faces, $f,g$, but if $f = g$, we still want $v$ to take charge from $f = g$ twice if $f$ is a $6^+$-face.
	
	Let $c^*$ be the resulting charge function after this redistribution.
	We will show that for all $v \in V(G)$, $c^*(v) \ge 0$ and for all $f \in F(G)$, $c^*(f) \ge 0$.

    \begin{description}
        \item[2-vertex] Let $v \in V(G)$ be a 2-vertex. Note $c(v) = -2$. Let $f,g$ be the faces incident to $v$. By Lemmas \ref{lem:no2v3f} and \ref{lem:no2v4f}, both $f$ and $g$ are $6^+$-faces. Therefore, $c^*(v) = c(v) + 2 = -2 + 2 = 0 \ge 0$.
        \item[3-vertex] Let $v \in V(G)$ be a 3-vertex. Note $c(v) = -1$. Let $f,g,h$ be the faces incident to $v$. By Lemmas \ref{lem:no3v_33f}, \ref{lem:no3v_44f}, and \ref{lem:no34f}, at most one of $f,g,h$ can be a $4^-$-face, and therefore the remaining two are $6^+$-faces. Therefore, $c^*(v) = c(v) + 2 \cdot \frac{1}{2} = -1 + 1 = 0 \ge 0$.
        \item[4-vertex] Let $v \in V(G)$ be a 4-vertex. Then $c^*(v) = c(v) = 0 \ge 0$.
        \item[3-face] Let $f \in F(G)$ be a 3-face. Note $c(f) = -1$. Let $g,h,k$ be the faces incident to $f$. By Lemmas \ref{lem:no333f} and \ref{lem:no34f}, we know that none of $g,h,k$ can be 4-faces and at most one can be a 3-face. Therefore, at least two of $g,h,k$ are $6^+$-faces, and so $c^*(f) = c(f) + 2 \cdot \frac{1}{2} = -1 + 1 = 0 \ge 0$.
        \item[4-face] Let $f \in F(G)$ be a 4-face. Note $f$ gives no charge and takes no charge and so remains unchanged. Then $c^*(f) = c(f) = 0 \ge 0$.
        \item[$6^+$-face] Let $f \in F(G)$ be a $6^+$-face. Let $l$ be the length of the face. Note $c(f) = l - 4$.
	    We will use intermediate discharging to implement the above rules. First, we will give each edge of the face $\frac{1}{3}$ charge. After this step, the final charge left in $f$ is:
	    \begin{align*}
	        l - 4 - \frac{l}{3} &= \frac{2l}{3} - 4 \\
	    \end{align*}
	    This is nonnegative since $l \ge 6$.
	    We will now distribute the charges from the edges around $f$ to the 2-vertices, 3-vertices, and 3-faces incident to $f$ that according to the above rules should be pulling charge from $f$, and each of these receives the charge that they should receive. 
    \end{description}
	
	\begin{enumerate}[SubR1]
	    \item A 3-face takes $\frac{1}{3}$ charge from the shared edge with $f$. (See Figure \ref{fig:srule3f}) \label{srule:3f}
	    \item A 3-face adjacent to at least one 3-face additionally takes $\frac{1}{6}$ charge from the edge on the same side as the incident 3-face. (See Figure \ref{fig:srule3f3f}) \label{srule:3f3f}
	    \item A 3-vertex incident to a 3-face takes $\frac{1}{3}$ charge from the other edge incident to it, and it takes $\frac{1}{6}$ charge from the edge on the other side of the 3-face along $f$. (See Figure \ref{fig:srule3v3f}) \label{srule:3v3f}
	    \item A 3-vertex not incident to a 3-face takes $\frac{1}{4}$ charge from each edge incident to it. (See Figure \ref{fig:srule3v}) \label{srule:3v}
	    \item A 2-vertex takes $\frac{1}{3}$ charge from the two edges incident to it, and takes $\frac{1}{6}$ charge from the edges distance two away from it along $f$. (See Figure \ref{fig:srule2v}) \label{srule:2v}
	\end{enumerate}
	
	\begin{figure}[ht]
    \centering
    \begin{tikzpicture}[node distance=2cm, every node/.style={draw=black, circle}]
        \node (A) at (0*360/3-30:1.5) {};
        \node (B) at (1*360/3-30:1.5) {};
        \node (C) at (2*360/3-30:1.5) {};
        \node[draw=none] (AC) at ($(A)!0.5!(C)$) {};
        \node[draw=none] (M) at ($(AC) + (0,1)$) {};
        \node[draw=none] at ($(AC) + (0,-0.5)$) {$f$};
        \draw (A) -- (B) -- (C) -- (A);
        \draw[-latex] (AC) -- (M) node[midway, right, draw=none] {$\frac{1}{3}$};
    \end{tikzpicture}
    \caption{A 3-face takes $\frac{1}{3}$ charge from its shared edge with $f$}
    \label{fig:srule3f}
\end{figure}
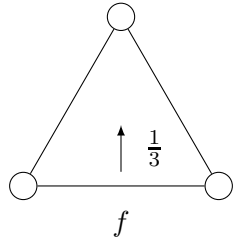

\begin{figure}[ht]
    \centering
    \begin{tikzpicture}[node distance=2cm, every node/.style={draw=black, circle}]
        \node (A) at (0*360/3-30:1.5) {};
        \node (B) at (1*360/3-30:1.5) {};
        \node (C) at (2*360/3-30:1.5) {};
        \node (D) at (360/6-30:3) {};
        \node[draw=none] (E) at ($2*(A) - (C)$) {$4$};
        \node[draw=none] (AC) at ($(A)!0.5!(C)$) {};
        \node[draw=none] (AE) at ($(A)!0.5!(E)$) {};
        \node[draw=none] (M) at ($(AC) + (0,1)$) {};
        \node[draw=none] (N) at ($(AC) + (0.1,1)$) {};
        \node[draw=none] at ($(AC) + (0,-0.5)$) {$f$};
        \draw (A) -- (B) -- (C) -- (A);
        \draw (B) -- (D) -- (A);
        \draw (A) -- (E);
        \draw[-latex] (AC) -- (M) node[midway, left, draw=none] {$\frac{1}{3}$};
        \draw[-latex] (AE) edge[bend left=60] node[midway, below, draw=none] {$\frac{1}{6}$} (N);
    \end{tikzpicture}
    \caption{A 3-face incident to at least one 3-face additionally takes $\frac{1}{6}$ charge from the edge on the same side as the incident 3-face.}
    \label{fig:srule3f3f}
\end{figure}
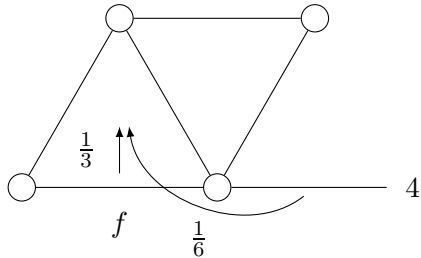

\begin{figure}[ht]
    \centering
    \begin{tikzpicture}[node distance=2cm, every node/.style={draw=black, circle}]
        \node (A) at (0*360/3-30:1.5) {};
        \node (B) at (1*360/3-30:1.5) {};
        \node[draw=none] (C) at (2*360/3-30:1.5) {$3$};
        \node[draw=none] (D) at ($2*(C) - (A)$) {$4$};
        \node[draw=none] (E) at ($2*(A) - (C)$) {$4$};
        \node[draw=none] (CD) at ($(C)!0.5!(D)$) {};
        \node[draw=none] (AE) at ($(A)!0.5!(E)$) {};
        \node[draw=none] at ($(A)!0.5!(C) + (0,-0.5)$) {$f$};
        \draw (A) -- (B) -- (C) -- (A);
        \draw (D) -- (C);
        \draw (A) -- (E);
        \draw[-latex] (CD) edge[bend left=60] node[midway, above, draw=none] {$\frac{1}{3}$} (C);
        \draw[-latex] (AE) edge[bend left=60] node[midway, below, draw=none] {$\frac{1}{6}$} (C);
    \end{tikzpicture}
    \caption{A 3-vertex incident to a 3-face takes $\frac{1}{3}$ charge from the other edge incident to it, and it takes $\frac{1}{6}$ charge from the edge on the other side of the 3-face along $f$}
    \label{fig:srule3v3f}
\end{figure}
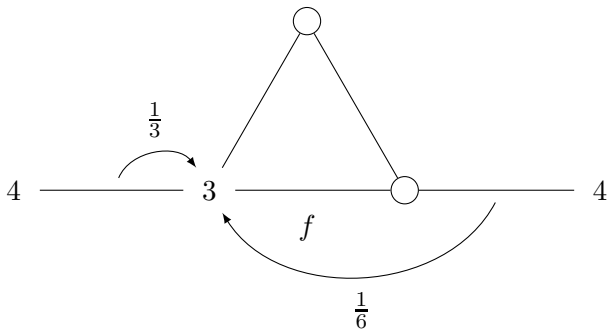

\begin{figure}[ht]
    \centering
    \begin{tikzpicture}[node distance=2cm, every node/.style={draw=black, circle}]
        \node[draw=none] (B) at (1.5,0) {$4$};
        \node[draw=none] (C) at (0,0) {$3$};
        \node[draw=none] (D) at (-1.5,0) {$4$};
        \node[draw=none] (BC) at ($(B)!0.5!(C)$) {};
        \node[draw=none] (CD) at ($(C)!0.5!(D)$) {};
        \node[draw=none] at ($(C) + (0,-1)$) {$f$};
        \draw (B) -- (C) -- (D);
        \draw[-latex] (BC) edge[bend left=60] node[midway, below, draw=none] {$\frac{1}{4}$} (C);
        \draw[-latex] (CD) edge[bend left=60] node[midway, above, draw=none] {$\frac{1}{4}$} (C);
    \end{tikzpicture}
    \caption{A 3-vertex not incident to a 3-face takes $\frac{1}{4}$ charge from each edge incident to it.}
    \label{fig:srule3v}
\end{figure}
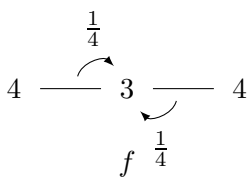

\begin{figure}[ht]
    \centering
    \begin{tikzpicture}[node distance=2cm, every node/.style={draw=black, circle}]
        \node[draw=none] (A) at (3,0) {$4$};
        \node[draw=none] (B) at (1.5,0) {$4$};
        \node[draw=none] (C) at (0,0) {$2$};
        \node[draw=none] (D) at (-1.5,0) {$4$};
        \node[draw=none] (E) at (-3,0) {$4$};
        \node[draw=none] (AB) at ($(A)!0.5!(B)$) {};
        \node[draw=none] (BC) at ($(B)!0.5!(C)$) {};
        \node[draw=none] (CD) at ($(C)!0.5!(D)$) {};
        \node[draw=none] (DE) at ($(D)!0.5!(E)$) {};
        \node[draw=none] at ($(C) + (0,-1)$) {$f$};
        \draw (A) -- (B) -- (C) -- (D) -- (E);
        \draw[-latex] (BC) edge[bend left=60] node[midway, below, draw=none] {$\frac{1}{3}$} (C);
        \draw[-latex] (CD) edge[bend left=60] node[midway, above, draw=none] {$\frac{1}{3}$} (C);
        \draw[-latex] (AB) edge[bend right=60] node[midway, above, draw=none] {$\frac{1}{6}$} (C);
        \draw[-latex] (DE) edge[bend right=60] node[midway, below, draw=none] {$\frac{1}{6}$} (C);
    \end{tikzpicture}
    \caption{A 2-vertex takes $\frac{1}{3}$ charge from the two edges incident to it, and takes $\frac{1}{6}$ charge from the edges distance two away from it along $f$.}
    \label{fig:srule2v}
\end{figure}
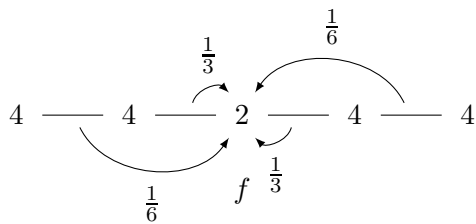
	
	We can see that these rules provide enough charge to satisfy the 2-vertices, 3-vertices, and 3-faces around $f$.
	
	Note that charges of $\frac{1}{3}$ and $\frac{1}{4}$ travel a ``short'' distance, from an edge to a vertex or face incident to it, while charges of $\frac{1}{6}$ travel a ``long'' distance, from an edge to a a vertex or face incident to an adjacent edge. Also note that since we are moving charge around a particular face $f$, there are only two directions to move it, restricting options significantly. 
	
	If an edge $e$ has two short distance charges leaving, or a short and a long distance charge leaving, or if it has two long distance charges leaving in the same direction, then it could end up with negative charge. We will show that this cannot happen.
	Lemmas \ref{lem:no2v3f}, \ref{lem:no22v}, \ref{lem:no23v}, \ref{lem:no33v}, and \ref{lem:no3v_33f} ensure that $e$ cannot have two short distance charges leaving. 
	Lemmas \ref{lem:no243v}, \ref{lem:no333f}, \ref{lem:no3v-3f3v}, and \ref{lem:no3v-m3f3f} show reducible the configurations which would occur if $e$ had a short distance charge and a long distance charge leaving, necessarily in different directions.
	Lemmas \ref{lem:no2v3f} and \ref{lem:no3v-m3f3f} show that the if $e$ has two long distance charges leaving in the same direction going to a vertex and a face, then there is a reducible configuration on that side. Since there cannot be any reducible configurations in $G$, we cannot have this happen.
	
	Therefore, we see that every edge ends up with nonnegative charge. Since every edge of the face $f$ ends up with nonnegative charge, and the face $f$, as shown above, ends up with nonnegative charge, the total is also nonnegative.
	
	We have shown that the final charge after redistribution is nonnegative for all vertices, edges, and faces of $G$. However, the initial charge was negative, and we only moved charge around. This is a contradiction. Therefore, the assumption that a minimal counterexample to our statement $G$ exists is wrong. And so the theorem is proved.
\end{proof}

\section{Future Work}

Future work will look towards extending Theorem \ref{thm:wegner4no5} to the more general statement for all planar graphs with max degree 4:

\begin{conjecture}
    \label{conj:wegner4}
    Let $G$ be a planar graph such that $\Delta(G) \le 4$. Then $\chi_\ell(G^2) \le 12$.
\end{conjecture}

\bibliographystyle{abbrv}
\bibliography{citations}

\begin{thebibliography}{10}

\bibitem{AppelHaken1977}
K.~Appel and W.~Haken.
\newblock {Every planar map is four colorable Part I: Discharging}.
\newblock {\em Illinois Journal of Mathematics}, 21(3):429--490, 1977.

\bibitem{AppelHaken1989}
K.~Appel and W.~Haken.
\newblock {\em Every Planar Map Is Four Colorable (Contemporary Mathematics)}.
\newblock American Mathematical Society, 1989.

\bibitem{AppelHakenKoch1977b}
K.~Appel, W.~Haken, and J.~Koch.
\newblock {Every planar map is four colorable Part II: Reducibility}.
\newblock {\em Illinois Journal of Mathematics}, 21(3):491--567, 1977.

\bibitem{Birkhoff1913}
G.~D. Birkhoff.
\newblock The reducibility of maps.
\newblock {\em American Journal of Mathematics}, 35(2):115, Apr 1913.

\bibitem{Borodin2013}
O.~V. Borodin.
\newblock {Colorings of plane graphs: A survey}.
\newblock {\em Discrete Mathematics}, 313(4):517--539, Feb 2013.

\bibitem{Borodin2002}
O.~V. Borodin, H.~J. Broersma, A.~Glebov, and J.~van~den Heuvel.
\newblock Stars and bunches in planar graphs. part {II}: General planar graphs
  and colourings.
\newblock Technical Report 0169-2690, University of Twente, Department of
  Applied Mathematics, 2002.

\bibitem{cranston_west_2017}
D.~W. Cranston and D.~B. West.
\newblock An introduction to the discharging method via graph coloring.
\newblock {\em Discrete Mathematics}, 340(4):766–793, 2017.

\bibitem{ErdosRubinTaylor1979}
P.~Erd\H{o}s, A.~L. Rubin, and H.~Taylor.
\newblock {Choosability in graphs}.
\newblock {\em Congressus Numerantium}, 26:125--127, 1979.

\bibitem{Heesch1969}
H.~Heesch.
\newblock {\em {Untersuchungen zum Vierfarbenproblem}}.
\newblock B.I.-Hochschulskripten, 810/810a/810b. Mannheim, Bibliographisches
  Institut, 1969.

\bibitem{Heesch1972}
H.~Heesch.
\newblock {Chromatic reduction of the triangulations $T_e$, $e = e_5 + e_7$}.
\newblock {\em Journal of Combinatorial Theory, Series B}, 13(1):46--55, Aug
  1972.

\bibitem{JendrolVoss2013}
S.~Jendrol' and H.~J. Voss.
\newblock {Light subgraphs of graphs embedded in the plane - A survey}.
\newblock {\em Discrete Mathematics}, 313(4):406--421, 2013.

\bibitem{Vizing1976}
V.~Vizing.
\newblock Vertex coloring of a graph with assigned colors (in russian).
\newblock {\em Metody Diskret. Analiz. (Novosibirsk)}, 29:3--10, 1976.

\bibitem{Wegner1977preprint}
G.~Wegner.
\newblock Graphs with given diameter and a colouring problem.
\newblock Preprint, 1977.

\bibitem{Zhu2018}
J.~Zhu and Y.~Bu.
\newblock Minimum 2-distance coloring of planar graphs and channel assignment.
\newblock {\em Journal of Combinatorial Optimization}, 36(1):55--64, Jul 2018.

\end{thebibliography}

\end{document}